\documentclass[12pt]{article}
\usepackage{amssymb}
\usepackage{amsmath}
\begin{document}

\title{A remark on a question of Beauville about lagrangian fibrations}

\author{Ekaterina Amerik\thanks{Universit\'e Paris-Sud, Campus d'Orsay,
B\^at. 425, 91405 Orsay, France; and National Research University Higher 
School of Economics, Department
of Mathematics, Vavilova 7, 117332, Moscow, Russia; 
Ekaterina.Amerik@math.u-psud.fr; the author is partially supported by
AG Laboratory, RF government grant ag. 11.G34.31.0023.}}

\date{}

\maketitle

\begin{abstract}

This note is a proof of the fact that a lagrangian torus on an irreducible 
hyperk\"ahler 
fourfold is always a fiber of an almost holomorphic lagrangian fibration.

\end{abstract}


\

Let $X$ be an irreducible holomorphic symplectic variety of dimension $2n$,
with a lagrangian torus $A\subset X$ (such a torus is always an abelian 
variety, even if $X$ is not projective: this follows from the fact that
the space $H^{2,0}(X)$ restricts to zero on $A$, see \cite{C}).
Beauville \cite{B} asked whether there always exists a lagrangian fibration $f$
(holomorphic or almost holomorphic), such that $A$ is a fiber of $f$. Recently, 
D. Greb, C. Lehn and S. Rollenske \cite{GLR} proved that it is always the
case when $X$ is non-projective, or even admits a non-projective deformation
preserving the lagrangian subtorus. It is not clear a priori whether all
irreducible holomorphic symplectic varieties admit such deformations.
The authors also have announced another paper in preparation, where 
Beauville's question is answered in the affirmative
in dimension four\footnote{Their paper appeared on arxiv.org:1110.2680 
a few days after a conversation between myself and one of the authors, and 
contains a simplified (with respect to 
the original approach of the authors) argument for the existence of a
meromorphic fibration, similar to the one 
given below. My primary reason for making this note public was to provide a
source which they could cite.}. The purpose of this short note is to provide an answer to
the
"almost holomorphic version" of Beauville's question in dimension four in a
very elementary way. 

\smallskip
 
Since the torus $A\subset X$ is lagrangian, its normal bundle $N_{A,X}$,
being isomorphic to its cotangent bundle, 
is trivial.
The deformation theory of lagrangian tori in holomorphic symplectic 
manifolds is well understood:
as it is explained in \cite{GLR}, by the results of Ran and Voisin 
the space $T$ of deformations of $A$ in $X$
is of dimension $n$, it is smooth at $A$, and any deformation $A_t$ of $A$ 
which
is smooth is itself a lagrangian torus. Moreover, if $Y$ denotes the
 universal family of the $A_t$'s, the projection
$q: Y\to X$ is unramified along each smooth $A_t$. In particular, 
one has a finite
number $d$ of $A_t$'s through a general point of $X$, and if through a certain
point $x\in X$ there is an infinite number of $A_t$ then almost all of
them are singular (in fact, if $Y$ denotes the universal family, the projection
$q: Y\to X$ is unramified along each smooth $A_t$). To show that there
is an almost holomorphic fibration with fibers $A_t$ is the same as to show that
$d=1$, or that a general $A_t$ does not intersect any other member $A_s$.

To say the same things using slightly different words, from the 
above-mentioned facts on the universal family 
map $q$ the
following is immediate.

\medskip

{\bf Fact:} {\it Locally in a neighbourhood $U$ of $A$, one has a lagrangian 
fibration by
its small deformations $A_t$.}

\medskip

One has to show that this local fibration gives rise to a global, possibly
meromorphic, fibration.

Let us list a few immediate consequences of this fact in a lemma.
We assume for simplicity that $X$ is projective.

\medskip

{\bf Lemma 1:} {\it a) A lagrangian torus does not intersect its
small deformations;

b) Whenever $A_t\cap A_s \neq \emptyset$ and $A_s$ is smooth, all 
irreducible
components of 
$A_t\cap A_s$ are positive-dimensional;

c) For a general $A_{s'}$ intersecting $A_t$,
the intersection $A_t\cap A_{s'}$ is equidimensional.}

\medskip

{\it Proof:} Whereas a) is immediate from the existence of a local
fibration, b) deserves a few words of explanation. Consider the local
fibration $f_s: U_s\rightarrow T_s$ in a neighbourhood of $A_s$ 
(here $T_s$ is a small disc around $s\in T$). The
intersection of $A_t$ with $U_s$ may consist of several local components.
The statement that $A_t\cap A_s$ cannot have a zero-dimensional component 
is clearly implied by the claim that no component of  $A_t\cap U_s$ 
dominates the base $T_s$ (indeed, such a component would have 
zero-dimensional general fiber under the map $f_s$, whereas for a 
non-dominating component all fibers are of strictly positive dimension). 
If this 
claim does not hold, we can choose a sufficiently general point 
$s'\in T_s$ outside of the
image of all non-dominating components of $A_t\cap U_s$. The resulting 
lagrangian
torus $A_{s'}$ has non-empty zero-dimensional intersection with $A_t$.
This is impossible since the 
self-intersection number $[A]\cdot [A]=0$.  

The proof of c) is similar to b): it suffices to choose $s'\in T_s$
sufficiently general in a "maximal" image of a component of $A_t\cap U_s$.
If this image is of dimension $c$, the intersection $A_t\cap A_{s'}$
is purely $(n-c)$-dimensional.

\medskip

Consider now the case $dim(X)=4$, then the intersection $A_t\cap A_s$ is a 
curve
if nonempty.
Fix a general lagrangian torus $A_t$. Since the family of its deformations  
$A_s, s\in T$ 
induces a local fibration in the neighbourhood of each smooth member, the intersections $A_t\cap A_s$ induce
a fibration or several fibrations on $A_t$. Indeed, since most of the 
$A_s$ 
intersecting $A_t$ are smooth, this is a
family of cycles on $A_t$ not intersecting their neighbours, as are the $A_s$ 
on $X$ itself; since these cycles are divisors on a surface, those are (possibly 
reducible) curves whose square is equal to zero, that is, our $A_t$ is
fibered in elliptic curves and our intersections are unions of fibers.

Of course, the $A_s$ intersecting $A_t$ do not necessarily form an
irreducible family, so apriori there can be several (say, $k>1$) such fibrations 
on $A_t$. Nevertheless an easy linear algebra argument shows that there
is in fact only one:

\medskip

{\bf Proposition 2:} {\it One has  $k=1$.}

\medskip

{\it Proof:} Otherwise, we can find two more tori $A_u$ and $A_s$ 
through a general point $p$ of $A_t$, in such a way that the pairwise
intersections $C_{ts}$, $C_{tu}$ and $C_{us}$ are curves with distinct
tangents at $p$ and the intersection
of all three tori has $p$ as an isolated point. One deduces easily from
Sard's lemma that for general $t,s$, the intersection $A_s\cap A_t$ is
reduced, so we may assume that the tangent planes to $A_u$, $A_s$ and $A_t$
are distinct. But since the pairwise intersections of those planes in $T_pX$ 
are distinct lines, the planes only span a hyperplane $V\subset T_pX$.
Now the restriction $\sigma_V$ of the symplectic form $\sigma_p$, that is,
the value of $\sigma$ at $p$, to $V$
has one-dimensional kernel. Recall that $T_pA_u$, $T_pA_s$ and $T_pA_t$ are
$\sigma_p$-isotropic. Since $V=T_pA_t+T_pA_s$, one must have 
$T_pC_{ts}=Ker(\sigma_V)$. But by the same reason the same holds for 
$T_pC_{tu}$, a contradiction.  

\medskip

Now we are ready to prove the announced result, which we formulate as
a theorem.

\medskip

{\bf Theorem 3:} {\it The answer to the almost holomorphic version of 
Beauville's 
question is positive in dimension 4, that is, a general $A_t$ does not
intersect any other $A_s$, so that $d=1$ and the family $A_t$ induces an
almost holomorphic fibration on $X$.}

\medskip

{\it Proof:} Suppose $d>1$, then we claim that the irreducible components 
of the intersections
$A_t\cap A_s$ rationally fiber $X$, that is, there is only one such 
curve through a general point $p\in X$. Indeed, suppose there are
two of them, $C_1$ and $C_2$. Then $C_1$ is a component of the intersection
$A_{1,1}\cap A_{1,2}$ and $C_2$ of $A_{2,1}\cap A_{2,2}$. Since $p$ is
general in $X$, we may suppose that all four of the above abelian surfaces 
are smooth. Each
pair of those four tori intersect at $p$, and we know that they should then 
intersect along a curve through $p$. Looking at the
possible configurations and keeping in mind that the intersection
of smooth surfaces is a disjoint union of elliptic curves, one is back 
to what is just ruled out in
Proposition 2.

On the other hand, it is proved in \cite{AC}, Proposition 3.4, that any
rational fibration $f: X\dasharrow B$ of an irreducible holomorphic 
symplectic variety $X$ has fibers of dimension at least $dim(X)/2$, 
provided that the general fiber is not of general type. Since the irreducible
components of $A_t\cap A_s$ are elliptic curves, they cannot
rationally fiber X, a contradiction.

\medskip

{\bf Remark 4:} This argument generalizes to higher-dimensional $X$ in the
case when the intersection of two general (intersecting) lagrangian tori 
$A_t\cap A_s$ is
of codimension one in each; so if $2n=dim(X)$, one wants $A_t\cap A_s$ to
be equal to $n-1$; when $n=2$ this is of course the only option unless 
$X$ is fibered by the $A_t$'s. More precisely, $A_t\cap A_s$ must be of
dimension $n-1$ for {\it each} general pair of intersecting lagrangian tori:
that is to say, the subset $W=\{(t,s)| A_t\cap A_s\neq \emptyset\}\subset T\times T$ 
can have
several irreducible components $W_i,\ i=1,\dots, l$, such that a generic
pair of tori in $W_i$ intersects in dimension $e_i$, and for our argument 
we need this dimension to be equal to $n-1$ for any $i=1,\dots, l$. 
As soon as this
condition is satisfied, the same reasoning as above applies to show the
existence of an almost holomorphic lagrangian fibration with fibers $A_t$.

\medskip

We hope to return to Beauville's question in a forthcoming joint work
with Fr\'ed\'eric Campana.

\end{document}